\documentclass[11pt]{amsart}
\usepackage{amssymb,amsfonts, amscd, amsmath}
\usepackage[all]{xy}
\CompileMatrices
\UseTips
\newtheorem{prop}{Proposition}
\newtheorem{lem}[prop]{Lemma}
\newtheorem{thm}[prop]{Theorem}

\newtheorem{defn}[prop]{Definition}
\newtheorem{ex}[prop]{Example}
\newtheorem{rem}[prop]{Remark}
\numberwithin{prop}{section}
\numberwithin{equation}{section}

\newcommand{\C}{{\Bbb C}}

\newcommand{\Q}{{\Bbb Q}}

\newcommand{\f}{{{\rm F}}}
\newcommand{\w}{{{\rm W}}}
\newcommand{\h}{{{\rm H}}}
\newcommand{\bs}{{B_s(\mathcal{A}^\bullet)}}

\renewcommand{\hom}{{\rm Hom}}
\newcommand{\ep}{\hfill \qedsymbol \newline }

\newcommand{\DR}{{\rm \footnotesize{DR}}}
\newcommand{\Dol}{\rm \footnotesize{Dol}}
\newcommand{\lra}{\longrightarrow}
\newcommand{\dbar}{{\bar{\partial}}}

\newcommand{\cinfty}{{\mathcal{C}^\infty}}

\begin{document}

\title[Unipotent flat bundles and Higgs bundles]
{Unipotent flat bundles and Higgs bundles over compact K\"ahler manifolds}

\author{Silke Lekaus}

\address{FB 6 - Mathematik, 
Universit\"at Essen,
45117 Essen,
Germany}

\email{silke.lekaus@uni-essen.de}


\maketitle

\begin{abstract}
We characterize those unipotent representations of the fundamental group $\pi_1(X,x)$ 
of a compact K\"ahler manifold $X$, which correspond
to a Higgs bundle whose underlying Higgs field is equal to zero.

The characterization is parallel to the one that R. Hain gave of those unipotent representations of $\pi_1(X,x)$ 
that can be realized as the monodromy of a flat connection on the holomorphically trivial vector bundle.

We see that Hain's result and ours follow from a careful study of Simpson's correspondence between Higgs bundles
and local systems.
\end{abstract}

\section{Introduction}

A Higgs bundle over a compact, complex K\"ahler manifold $X$ is a pair $(E,\theta)$, where $E$ is a holomorphic
vector bundle over $X$ and $\theta : E\to E\otimes \Omega^1_X$ is a holomorphic map
such that $\theta \wedge \theta =0$.
Higgs bundles play an important role in non-abelian Hodge theory.

In \cite{si}, C. Simpson proves that there is an equivalence of categories between the category of flat bundles
and the category of Higgs bundles that are extensions of stable Higgs bundles with vanishing Chern classes.

R. Hain investigated when a unipotent group homomorphism $\rho : \pi_1(X,x)$ $\to {\rm Gl}(n,\C)$,
$x\in X$, can be realized as the monodromy representation
of a flat connection on the holomorphically trivial bundle.

Since $\rho$ is assumed to be unipotent, its $\C$-linear extension to $\C\pi_1(X,x)$ factors through a power $J^{s+1}$
of the augmentation ideal (see Section \ref{unipotentsection}).
The quotient $\C\pi_1(X,x)/J^{s+1}$ carries a mixed Hodge structure as defined by Morgan (\cite{morgan})
and Hain (\cite{hain1}). 

Hain proved in \cite{hain0} that a unipotent $\rho$ can be realized as the monodromy of a flat connection
on the holomorphically trivial bundle, if and
only if its $\C$-linear extension factors through the ideal 
\[
I := J\cap \f^0(\frac{\C\pi_1(X,x)}{J^{s+1}}) + J^2\cap \f^{-1}(\frac{\C\pi_1(X,x)}{J^{s+1}}) 
 + \cdots + J^s\cap \f^{-(s-1)}(\frac{\C\pi_1(X,x)}{J^{s+1}}).
\]

Using Simpson's correspondence Hain's result translates into

\begin{thm} (Hain, \cite{hain0})  \label{theorem}
Let $(V,\nabla)$ be a flat bundle whose $\C$-linear extension $\rho$ of the monodromy representation 
factors through $J^{s+1}$.

Then $\rho$ factors through $I$ if and only if the underlying holomorphic bundle of any representative
of the isomorphism class of Higgs bundles corresponding to $(V,\nabla)$ is trivial.
\end{thm}

We define
\[
\bar{I} = J\cap \bar{\f}^0(\frac{\C\pi_1(X,x)}{J^{s+1}}) + J^2\cap \bar{\f}^{-1}(\frac{\C\pi_1(X,x)}{J^{s+1}})
 + \cdots + J^s\cap \bar{\f}^{-(s-1)}(\frac{\C\pi_1(X,x)}{J^{s+1}}).
\]

The main result of this article is 
\begin{thm}  \label{bartheorem}
Let $(V,\nabla)$ be a flat bundle whose $\C$-linear extension $\rho$ of the monodromy representation 
factors through $J^{s+1}$. 

Then $\rho$ factors through $\bar{I}$ if and only if any representative 
of the isomorphism class of Higgs bundles corresponding to $(V,\nabla)$ has Higgs field equal to zero.
\end{thm}

The Sections \ref{bundlesection}, \ref{unipotentsection} and \ref{MHSsection} serve to provide the reader with
background material about the correspondence proved by Simpson, unipotent representations and the mixed Hodge 
structure on 
$\C\pi_1(X,x)/J^s$ and to fix the notations used in later paragraphs. 
There are no new results in these sections and they can be skipped by the reader familiar with these topics.

In Section \ref{J2case}, we shortly describe the special case of bundles whose monodromy factors through $J^2$,
which illustrates the method of the proofs in the general case, which is done in Section \ref{Jncase}.

\section{Correspondence between flat bundles and Higgs bundles} \label{bundlesection}

Let $X$ be a compact complex K\"ahler manifold. 
A Higgs bundle $(E,\theta)$ over $X$ as defined in the introduction, 
has an equivalent $\mathcal{C}^\infty$-description (see \cite{si}):
It can be considered as a
$\mathcal{C}^\infty$-bundle with a first order operator $D''=\dbar+\theta$,
where $\dbar$ defines the holomorphic structure of $E$ by taking sections
of $E$ to $(0,1)$-forms with coefficients in $E$ and by annihilating
holomorphic sections. It holds that $D''$ is integrable, i.e.
$(D'')^2=0$, because $\dbar^2=0$,
$\theta \wedge \theta =0$, and $\dbar(\theta)=0$, as $\theta$ is
holomorphic. Further it fulfills the Leibniz rule
\[D''(fe)=\dbar(f)\wedge e+fD''(e)\]
for all sections $f$ of $\mathcal{A}^0$ and all sections $e$ of $E$. 

Conversely, let $D'' =\dbar +\theta$ be an operator on a
$\mathcal{C}^\infty$-vector bundle $E$,
so that $\dbar$ takes sections of $E$ to $(0,1)$-forms
with coefficients in $E$ and $\theta$ takes sections of $E$ to
$(1,0)$-forms with coefficients in $E$. Such an operator defines a Higgs
structure on $E$, if and only if it is integrable and
fulfills the Leibniz rule. These conditions imply that
$\dbar^2=0$, $\theta\wedge\theta=0$, and that $\theta$ is holomorphic.

The categories of flat bundles and Higgs bundles over $X$ form so called differential graded
categories, if we define an appropriate Hom-complex (see \cite{si}):

\begin{ex} \label{catexamples}
{\rm (1) By
$\mathcal{C}_\DR$ we denote the differential graded category whose objects
are all flat bundles over $X$ and which has the Hom-complex
\[\mbox{Hom}^\bullet(U,V)=
(\Gamma(X,\mathcal{A}^\bullet(\mbox{Hom}(U,V))),D),\]
where the composition of homomorphisms is given by the wedging of forms.

(2) By $\mathcal{C}_\DR^s$ we denote the full subcategory of $\mathcal{C}_\DR$
consisting of semisimple objects.

(3) By $\mathcal{C}_{\Dol}$, we denote the differential graded category whose
objects are those Higgs bundles on $X$ which are extensions of stable Higgs
bundles with vanishing Chern classes.
Its Hom-complex is
\[\mbox{Hom}^\bullet(U,V)=
(\Gamma(X,\mathcal{A}^\bullet(\mbox{Hom}(U,V))),D''),\]
where the composition of morphisms is obtained by the wedging of forms.

(4) By $\mathcal{C}_{\Dol}^s$ we denote the full subcategory consisting of semisimple
objects, i.e. polystable Higgs bundles.

(5) By $\mathcal{C}_{\rm harm}^s$ we denote the category of all harmonic bundles (i.e.
those bundles which carry both a flat and a Higgs structure) with
Hom-complex
\[\mbox{Hom}^\bullet(U,V)=({\rm ker}(D')\subset
\Gamma(X,\mathcal{A}^\bullet(\mbox{Hom}(U,V))),D''),\]
where $D$ is the flat operator and $D''=\dbar + \theta$ the Higgs operator on Hom$(U,V)$,
and $D':=D-D''$.} \ep
\end{ex}

If $\mathcal{C}$ is a differential graded category with differential $d$ on the Hom-complex, 
we can construct its completion $\hat{\mathcal{C}}$ (see \cite{gm}, \cite{si}):

Let $\bar{\mathcal{C}}$ be the differential graded category,
whose objects are pairs $(U,\eta)$ with an object $U$ of $\mathcal{C}$
and an endomorphism $\eta\in\mbox{Hom}^1(U,U)$ of degree one, satisfying
$d(\eta) + \eta^2 =0$. If $(U,\eta)$ and $(V,\xi)$ are two objects of
$\bar{\mathcal{C}}$, we define the differential graded algebra 
Hom$((U,\eta),(V,\xi))$ to be Hom$(U,V)$
with the same grading, the new differential
$\hat{d}(f) = d(f) + \xi f - (-1)^{\deg(f)}f\eta.$

The category $\mathcal{C}$ can be embedded into $\bar{\mathcal{C}}$ by
mapping an object $U$ to $(U,0)$. 

We define the {\it completion} $\hat{\mathcal{C}}$ of $\mathcal{C}$
to be the full subcategory
of  $\bar{\mathcal{C}}$, whose objects are successive
extensions of objects of $\mathcal{C}$.

There is the following lemma:

\begin{lem} (Simpson, \cite{si})
There are equivalences of differential graded categories
\[\hat{\mathcal{C}}_{\DR}^s \cong \mathcal{C}_{\DR} \, \mbox{ and } \,
\hat{\mathcal{C}}_{\Dol}^s \cong \mathcal{C}_{\Dol}.\]
\end{lem}

\begin{rem}
{\rm In the following we will need a concrete description of the functors
\[ \hat{\mathcal{C}}_{\DR}^s \to \mathcal{C}_{\DR} \, \mbox{ and } \,
\hat{\mathcal{C}}_{\Dol}^s \to \mathcal{C}_{\Dol},\]
which can be found in \cite{si}, Section 3:\\
Let a pair $((U,D),\eta)$ be an object of $\hat{\mathcal{C}}_{\DR}^s$, i.e. 
$(U,D)$ is a semisimple flat bundle and $\eta$ is 1-form with values in
End$(U)$ satisfying $D(\eta) + \eta\wedge\eta=0$. The corresponding (in general not 
semisimple) object of $\mathcal{C}_{\DR}$ is $(U,D + \eta)$.\\
Analogously, an object $((U,D''),\eta)$ of $\hat{\mathcal{C}}_{\Dol}^s$
is mapped to $(U, D'' + \eta)$. }
\end{rem}

If $(\mathcal{C},d)$ is a differential graded category, we denote by E$^0(\mathcal{C})$ the additive
category whose objects are the same as those of $\mathcal{C}$, but whose morphisms from $U$ to $V$
are just Ext$^0(U,V)={\rm Ker }\, d:{\rm Hom}^0(U,V)\to {\rm Hom}^1(U,V)$.

\begin{lem} (Simpson, \cite{si}) \label{Simpsoncorrespondence} 
There is an equivalence of categories between the category {\rm E}$^0(\mathcal{C}_{\DR})$ of flat bundles and
the category {\rm E}$^0(\mathcal{C}_{\Dol})$ of Higgs bundles which are extensions of stable Higgs bundles with 
vanishing Chern classes. Furthermore these are equivalent to the category
{\rm E}$^0(\hat{\mathcal{C}}^s_{\rm harm})$.
\end{lem}

\begin{rem} \label{Corlette}
{\rm The equivalence E$^0(\mathcal{C}_{\DR})\,\cong\,{\rm E}^0(\mathcal{C}_{\Dol})$ can be described by using
the equivalence to E$^0(\hat{\mathcal{C}}^s_{\rm harm})$ (as suggested by Corlette, see \cite{si}):

An object of E$^0(\hat{\mathcal{C}}^s_{\rm harm})$ is a pair $(U,\eta)$ with a harmonic bundle
$U=(U,D,D',D'')$  and $\eta\in\Gamma(X,\mathcal{A}^1({\rm End}(U))$ such that $D'(\eta)=0$, and
$D''(\eta) + \eta\wedge\eta=0$. A pair $(U,\eta)$ is mapped to the Higgs bundle $(U,D''+\eta)$ and
to the flat bundle $(U,D+\eta)$. Since these functors are equivalences of categories, for any flat
bundle and for any Higgs bundle, satisfying the conditions of the lemma, there is a unique pair
$(U,\eta)$ which maps to it.

Note: Since a harmonic bundle is semisimple as a flat bundle (see \cite{si}, Theorem 1),
a unipotent harmonic bundle $U$ is necessarily isomorphic to a direct sum of trivial bundles
$(\mathcal{A}^0_X,d)^{\oplus r}$.
Thus there exists a global flat frame for $U$, and with respect to this frame 
$D=d$, $D'=\partial$, and $D''=\dbar$ are the
usual differential operators operating coordinatewise. }
\end{rem}

\section{Unipotent representations and unipotent flat bundles} \label{unipotentsection}

Let $G$ be a group.
The homomorphism of $G$ to the trivial group induces the algebra
homomorphism
\[\begin{array}{cclcl}
\epsilon & : & \C G & \to & \C, \\
 & & \sum \alpha_g g & \to & \sum \alpha_g.
\end{array}\]

The kernel
\[J:=\mbox{Ker }(\epsilon)\]
is called the {\it augmentation ideal}.
It is spanned by $\{g-1 \, | \, g\in G\}$.  \\

Let $V$ be a finite-dimensional $\C$-vector space and
 $\rho : G\to \mbox{Aut}(V)$ a representation. It can be extended by linearity
to an algebra homomorphism
\[\bar{\rho}: \C G \to \mbox{End}(V).\]

\begin{prop} \label{unipotentprop}
The following are equivalent:
\begin{enumerate}
\item $\rho$ is unipotent.
\item $\rho$ induces a homomorphism
\[\bar{\rho}: \C G / J^{s+1}\to \mbox{End}(V)\]
for some $0\le s \le \dim V$.
\end{enumerate}
\end{prop}

\begin{proof}
This follows from Kolchin's theorem (see e.g. \cite{se}, Part I, Chapter V), saying that there is a basis for
$V$, in which all elements of $\rho(G)$ are simultaneously represented
by upper-triangular matrices with 1's on the diagonal.
\end{proof}

In the following $X$ denotes again a compact complex K\"ahler manifold.

\begin{defn}
A flat bundle $(V,\nabla)$ over $X$ is called unipotent, if its
monodromy representation $\rho$ is unipotent.
\end{defn}

The unipotent flat bundles are exactly those, which are successive
extensions of $(\mathcal{A}^0_X,d)$ by itself, as follows from
the following proposition:
\begin{prop} \label{RHprop}
There is an equivalence of categories between the category of unipotent
representations of $\pi_1(X,x)$ and the category of flat bundles, which are
successive extensions of $(\mathcal{A}^0_X,d)$.\\
The correspondence is the usual Riemann-Hilbert correspondence.
\end{prop}
\begin{proof}
This holds since the Riemann-Hilbert correspondence respects exact
sequences.
\end{proof}

\section{The reduced bar construction and the mixed Hodge structure on $\h^0(\bs)$} \label{MHSsection}

We explain the reduced bar construction following \cite{hain2}.

Let $A^+:=\bigoplus_{i=1}^\infty \Gamma(X,\mathcal{A}^i)$. For each $s\ge 0$ we define
\[T^{-s,t}:=[\otimes^s A^+]^t ,\]
where $t$ is the sum of the degrees of the elements in $A^+$. We denote an element $\omega_1\otimes \cdots \otimes \omega_s$
by $[\omega_1 | \ldots | \omega_s]$ and give it the total degree $s-t$.

On every $T^{-s,t}$ we define two differentials, namely the {\bf combinatorial  differential}
$d_C:T^{-s,t} \to T^{-s+1,t}$
given by
\[d_C([\omega_1|\ldots | \omega_s]):= \sum_{i=1}^{s-1}(-1)^{i+1}[J\omega_1|\dots |J\omega_{i-1}|J\omega_i\wedge \omega_{i+1}| \omega_{i+2}|\ldots | \omega_s],\]
and the {\bf internal differential}
$d_I:T^{-s,t} \to T^{-s,t+1}$
given by
\[d_I([\omega_1|\ldots | \omega_s]):= \sum_{i=1}^{s}(-1)^{i}[J\omega_1|\dots |J\omega_{i-1}|d\omega_i| \omega_{i+1}|\ldots | \omega_s],\]
where $J:A^+ \to A^+$ denotes the involution $J(\omega):=(-1)^{\deg \omega}\cdot \omega$.

It holds that 
\[d_I^2=d_C^2=0, d_I d_C+d_C d_I=0,\]
such that $(T^{\bullet\bullet}, d_I, d_C)$
is a double complex. 

We denote the total complex by $T(\mathcal{A}^\bullet)$.

We define $R$ to be the graded subspace of $T(\mathcal{A}^\bullet)$ spanned by the elements of the form
\begin{enumerate}
\item[(a)] $[df|\omega_1|\ldots|\omega_s] +  f(x)[\omega_1|\ldots|\omega_s] -  [f\omega_1|\ldots|\omega_s]$,
\item[(b)] $[\omega_1|\ldots|\omega_{i-1}|df|\omega_i|\ldots|\omega_s] +  [\omega_1|\ldots|f\omega_{i-1}|\omega_i|\ldots|\omega_s]$
\item[] \hspace{0.4cm} $ - [\omega_1|\ldots|\omega_{i-1}|f\omega_i|\ldots|\omega_s]$,
\item[(c)] $[\omega_1|\ldots|\omega_s|df] + [\omega_1|\ldots|f\omega_s] - [\omega_1|\ldots|\omega_s]f(x).$
\end{enumerate}

The {\bf reduced bar construction} $\bar{B}(\mathcal{A}^\bullet)$ is defined to be the quotient 
$T(\mathcal{A}^\bullet)/R$.

Defining
\[\mathcal{T}^{-s}=\oplus_{u\le s}T^{-u,v},\]
for each $s\ge 0$, we obtain the filtration 
\[ T(\mathcal{A}^\bullet)\supseteq \cdots \supseteq \mathcal{T}^{-2}\supseteq \mathcal{T}^{-1}\supseteq \mathcal{T}^0=\C\supseteq 0 \]
of $T(\mathcal{A}^\bullet)$, which induces the {\bf bar filtration} 
\[\bar{B}(\mathcal{A}^\bullet) \supseteq \cdots \supseteq \bar{B}_{2}\supseteq \bar{B}_{1}\supseteq \bar{B}_0=\C\supseteq 0 \]
of $\bar{B}(\mathcal{A}^\bullet)$.

We denote by $\h^0(\bs):=\C\oplus\h^0(\bar{B}_s(\mathcal{A}^\bullet))$ the set of classes of iterated integrals of at most 
length $s$ which are homotopy 
functionals, when restricted to loops at $x$. 

An iterated integral is a sum
\begin{equation}  \label{itint}
c + \sum_{k=1}^{n_1}\int \omega^{(1)}_k + \sum_{k=1}^{n_2}\int\omega^{(2)}_{k1}\omega^{(2)}_{k2}
+ \cdots + \sum_{k=1}^{n_s}\int\omega^{(s)}_{k1}\omega^{(s)}_{k2}\cdots \omega^{(s)}_{ks}
\end{equation}
with $c\in \C$ and $\omega_{ki}^{(j)}\in\Gamma(X,\mathcal{A}^1_X)$.

It is homotopy invariant, when restricted to loops based at $x$, if and only if 
\[(d_C+d_I)(\sum_{k=1}^{n_1}[\omega^{(1)}_k] + \sum_{k=1}^{n_2}[\omega^{(2)}_{k1}|\omega^{(2)}_{k2}]
+ \cdots + \sum_{k=1}^{n_s}[\omega^{(s)}_{k1}|\omega^{(s)}_{k2}|\ldots |\omega^{(s)}_{ks}]) \in R.\]

We consider such integrals to be functions on $\C\pi_1(X,x)$. As such they are not uniquely determined, since
two iterated integrals $I$ and $I'$ can give the same value on all loops based at $x$. This is the case,
if and only if their difference $I-I'$ lies in the $\C$-subvector space generated by 
the following three types of integrals (see e.g. \cite{hain1}, Proposition (1.3)):
\begin{enumerate}
\item[(a)] $\int df \omega_1\cdots\omega_r - \int (f\omega_1)\omega_2\cdots \omega_r +f(x)\int\omega_1\cdots\omega_r$,
\item[(b)] $\int\omega_1\cdots\omega_{i-1}df\omega_i\cdots\omega_r-\int\omega_1\cdots\omega_{i-1}(f\omega_i)\omega_{i+1}
\cdots\omega_r$\\
\hspace*{0.4cm}$ + \int \omega_1\cdots(f\omega_{i-1})\omega_i\cdots \omega_r$,
\item[(c)] $\int\omega_1\cdots \omega_r df - f(x)\int\omega_1\cdots\omega_r+\int\omega_1\cdots\omega_{r-1}
(f\omega_r)$
\end{enumerate}
with $\omega_1,\dots,\omega_r\in\Gamma(X,\mathcal{A}^1_X)$ and $f\in\Gamma(X,\mathcal{A}^0_X)$.\\
This is equivalent to saying that the involved 1-forms define an element of $R$.\\

There is a mixed Hodge structure on $\h^0(\bs)$ given by the following filtrations (for details see \cite{hain1}):
With $\f^p B(\mathcal{A}^\bullet)$ denoting the set of $\C$-linear combinations of iterated integrals,
such that each summand contains at least $p$ $dz's$, we define
\begin{eqnarray}
\f^p\h^0(B_s(\mathcal{A}^\bullet)) & :=  & B_s\h^0(\f^p B(\mathcal{A}^\bullet)),\label{Hfiltration}\\
\w_l\h^0( B_s(\mathcal{A}^\bullet)) & := & \h^0( B_l(X)),\label{Wfiltration}
\end{eqnarray}
i.e. the elements of $\f^p\h^0( B_s(\mathcal{A}^\bullet))$ (or $\bar{\f}^p\h^0( B_s(\mathcal{A}^\bullet))$) are those which are represented by
an iterated integral in $\f^p B(\mathcal{A}^\bullet)$ (or in $\bar{\f}^p B(\mathcal{A}^\bullet)$).
The elements of $\w_l\h^0( B_s(\mathcal{A}^\bullet))$ are those having a representative of length at most $l$, which means that
every summand of the representating iterated integral involves at most $l$ 1-forms.

\begin{thm} (Chen)
There is an isomorphism
\[\h^0(\bs) \lra \hom(\frac{\C\pi_1(X,x)}{J^{s+1}},\C),\]
given by integrating (classes of) iterated integrals over elements of $\C\pi_1(X,x)$.
\end{thm}
This is a well-defined map, since the integration of an iterated integral of length at most s
over an element of $J^{s+1}$ gives zero.

This isomorphism induces a dual mixed Hodge structure on $\C\pi_1(X,x)/J^{s+1}$, via
\begin{eqnarray} \label{Hodgefiltration}
\f^p\hom_\C 
 & = & \{\Phi\, |\, \Phi(\f^{1-p}(\C\pi_1(X,x)/J^{s+1}))=0 \}
\end{eqnarray} 
and
\begin{eqnarray}  \label{weightfiltration}
\w_l\hom_\Q 
 & = & \{\Phi\, |\, \Phi(\w_{-1-l}(\C\pi_1(X,x)/J^{s+1}))=0 \}.
\end{eqnarray}

\section{Illustration: The $J^2$-case}  \label{J2case}

Let $(V,\nabla)$ be a flat bundle of rank $r$ with monodromy representation $\rho$, whose $\C$-linear
extension, again denoted by $\rho$, factors through $J^2$.
There is an isomorphism class of Higgs bundles corresponding to $(V,\nabla)$
via Simpson's correspondence Lemma \ref{Simpsoncorrespondence}. 

Since $\rho$ factors through $J^2$, $(V,\nabla)$ is an extension of some $(\mathcal{A}^0_X,d)^{\oplus r_1}$
by $(\mathcal{A}^0_X,d)^{\oplus r_2}$ with $r_1+r_2=r$.

By Remark \ref{Corlette}, $(V,\nabla)$ is of the form $(V,D+\eta)$, where $(V,D)$ is a harmonic
bundle (hence semi-simple) and $\eta$ is in $\Gamma(X,\mathcal{A}^1_X({\rm End}(V)))$, satisfying 
\begin{equation}
D' \eta=0 \, \mbox{ and } \, D'' \eta + \eta\wedge\eta=0.
\end{equation}
Since $(V,D)$ is unipotent, it must hold that $(V,D)\cong (\mathcal{A}^0_X, d)^{\oplus r}$. 

Hence there is a global flat $\cinfty$-frame $f=(f_1,\ldots , f_r)$ for $(V,D)$,
with respect to which $D=d$, 
hence $D'=\partial$ and $D''=\dbar$, applied coordinatewise.

With respect to the same frame the connection matrix $A$ of $\nabla=d + A$ has the form
\begin{equation} \label{omega}
A = \left( \begin{array}{cc}
0 & \omega \\
0 & 0
\end{array}\right),\end{equation}
where $\omega=(\omega_{ij})$ is a $r_1\times r_2$-matrix of 1-forms.
The integrability condition implies that all 1-forms $\omega_{ij}$ are
closed.

The operator $\eta=\nabla - D$ must be given by the connection matrix $A$.

The operators $D$, $D'$, and $D''$, operating on $\eta$, are the induced ones on End$(V)$.
Therefore the application of the operators $D'$ and $D''$ on $\eta$ is also the coordinatewise application of 
$\partial$ and $\dbar$,
and the condition (5.1) implies that $\partial\omega_{ij}=0$ and $\dbar\omega_{ij}=0$, i.e.
the $\omega_{ij}$ are harmonic 1-forms.

Now $(V, D,D',D'',\eta)$ is the bundle in the category
$\hat{\mathcal{C}}^s_{\rm harm}$, defined in Example \ref{catexamples}, corresponding to $(V,\nabla)$.

The corresponding Higgs bundle, via Simpson's correspondence \ref{Simpsoncorrespondence}, is
\begin{equation}
(V,D''+\eta)=
(\bigoplus_{i=1}^r f_i \mathcal{A}^0_X, \dbar + \left( \begin{array}{cc}
0 & \omega\\
0 & 0
\end{array}\right) ).\end{equation}

The following lemma is a special case of \ref{typelemma}:

\begin{lem} \label{easytypelemma}
\begin{enumerate}
\item The underlying holomorphic bundle of $(V,D''+\eta)$ is trivial if and only if $\eta$ is of type (1,0),
i.e. if and only if all the forms $\omega=(\omega_{ij})$ in the formula (\ref{omega}) are of type (1,0). 
\item
The Higgs field of $(V,D''+\eta)$ is equal to zero if and only if $\eta$ is of type (0,1), i.e.
if and only if all the forms $\omega=(\omega_{ij})$ in the formula (\ref{omega}) are of type (0,1). 
\end{enumerate}
\end{lem}

With the notation above the monodromy representation induces the well-defined homomorphism
\[\begin{array}{ccccl}
\bar{\rho} & : & \frac{\C\pi_1(X,x)}{J^2} & \longrightarrow & {\rm M}(r,\C)\\
           &   &  [c] & \lra & {\rm E}_r + \left( \begin{array}{cc}
                                                 0 & \int_c \omega_{ij}\\
                                                 0 & 0
                                                \end{array}\right).
\end{array}\]

Since in this case the Hodge filtration on $\C\pi_1(X,x)/J^2$ is
\[\f^i(\C\pi_1(X,x)/J^2) = \left\{
\begin{array}{lcl}
 0, & & i=1,\\
\C\oplus\{c\in J/J^2 \, | \, \int_c\omega=0 \, \, \forall \mbox{ hol. }\omega\},& & i=0,\\
\C\pi_1(X,x)/J^2, & & i=-1, 
\end{array}\right. \]
and the ideals $I$ and $\bar{I}$ as defined in the introduction are
\begin{equation}
I =  \{c\in J/J^2 \, | \, \int_c\omega = 0 \mbox{ for all holomorphic }\omega\},
\end{equation}
and  
\begin{equation}
\bar{I} =  \{c\in J/J^2 \, | \, \int_c\omega = 0 \mbox{ for all anti-holomorphic }\omega\},
\end{equation} 
the if-parts of Theorems \ref{theorem} and \ref{bartheorem} are easily implied by
Lemma \ref{easytypelemma}. 

Conversely, the factorization through $I$ (or $\bar{I}$) implies 
that $\int\omega_{ij}\in \f^1\h^0$ (or in $\bar{\f}^1\h^0$). 
By the duality of the mixed Hodge structures we first know that there are $\cinfty$-functions $\psi_{ij}$ 
such that $\omega_{ij}-d\psi_{ij}$
are holomorphic (or anti-holomorphic), since there must be a representative of $\int\omega_{ij}$ in 
$\h^0(B_1(\mathcal{A}^\bullet))$ with this property.
As the $\omega_{ij}$ are harmonic, it follows from the  
Hodge decomposition theorem that $d\psi_{ij}=0$, such that Lemma \ref{easytypelemma} implies the assertion.

\section{The general case}  \label{Jncase}

Let $(V,\nabla)$ be a flat bundle of rank $r$ with monodromy representation $\rho$, whose $\C$-linear
extension, which we denote again by $\rho$, factors through $J^{s+1}$. 
By Propositions \ref{unipotentprop} and \ref{RHprop}, 
such a bundle is unipotent, and equivalently
a successive extension of direct sums $(\mathcal{A}^0_X,d)^{\oplus r_i}$, $i=1,\ldots, s+1$, of trivial flat 
bundles with $r=\sum_{i=1}^{s+1} r_i$. 
 
By Remark \ref{Corlette}, we know that $(V,\nabla)$ is of the form $(V,D+\eta)$, where
$(V,D,D',D'')$ is a harmonic bundle and $\eta\in\Gamma(X,\mathcal{A}^1_X({\rm End}(V)))$, satisfying
\begin{equation} \label{etaconditions}
D' \eta=0 \, \mbox{ and } \, D'' \eta + \eta\wedge\eta=0.
\end{equation}
By \cite{si}, Theorem 1, $(V,D)$ is semisimple.
Since it is in addition unipotent, it must hold that $(V,D)\cong (\mathcal{A}^0_X, d)^{\oplus r}$.

Therefore there is a global $D$-flat frame $f=(f_1,\ldots , f_r)$ for $V$,
and with respect to this frame $D=d$, 
hence $D'=\partial$ and $D''=\dbar$, applied coordinatewise.

With respect to the same frame, $\nabla$ is of the form $d+A$ with a connection matrix
\begin{equation} \label{connectionmatrix}
A=\left(\begin{array}{ccccc}
0 & \omega_{12}& \omega_{13} & \cdots  & \omega_{1,s+1}\\
0 & 0 & \omega_{23} & \cdots & \omega_{2,s+1}\\
\vdots & & & & \vdots\\
0 & 0 & 0 & \cdots & \omega_{s,s+1}\\
0 & 0 & 0 & \cdots & 0
\end{array}\right),
\end{equation}
where the $\omega_{ij}$ are $r_i \times r_j$-matrices of 1-forms,
and the operator $\eta=\nabla - D$ is given by the connection matrix $A$.

The conditions (\ref{etaconditions}) become 
\begin{equation} \label{etaconditions2}
\partial A =0 \, \mbox{ and }\, \dbar A + A\wedge A=0,
\end{equation}
where $\partial$ and $\dbar$ are applied coordinatewise.

The corresponding Higgs bundle, via Simpson's correspondence Lemma \ref{Simpsoncorrespondence}, is
\begin{equation}
(V,D''+\eta)=
(\bigoplus_{i=1}^r f_i \mathcal{A}^0_X, \dbar + A).
\end{equation}

\begin{lem}  \label{typelemma}
\begin{enumerate}
\item The underlying holomorphic bundle of $(V,D''+\eta)$ is trivial if and only if $\eta$ is of type (1,0),
i.e. if and only if the matrix $A$ is of type (1,0). 
\item
The Higgs field of $(V,D''+\eta)$ is equal to zero if and only if $\eta$ is of type (0,1), i.e.
if and only if the matrix $A$ is  of type (0,1). 
\end{enumerate}
\end{lem}
\begin{proof}
The assertion (2) is trivial. Further it is clear that the underlying holomorphic bundle of $(V,D'')$ is 
trivial, if $\eta$ is of type (0,1).
Conversely, we denote the (0,1)-part of $\eta$ by $\eta^{0,1}$, and we
assume that there exist global sections $h_1,\ldots,h_r \in \Gamma(X,V)$ satisfying 
$(D''+\eta^{0,1})(h_i)=0$ for all $i$ and forming a global frame for $V$.
With respect to the frame $f$, fixed above, the $h_i=\sum_{l=1}^r \tilde{\phi}_{il}f_l$ have $\cinfty$-coordinates
$\tilde{\phi}_{il}$.
Writing
\[\phi_{ij}=\left( \begin{array}{c}
\tilde{\phi}_{i,r_1+\ldots+r_{j-1}+1} \\
\vdots\\
\tilde{\phi}_{i,r_1+\ldots+r_j}
\end{array}\right)
\]
they satisfy the differential equation

\begin{equation} \label{diffeq}
\left( \begin{array}{c}
\dbar \phi_{i1} \\
\dbar \phi_{i2} \\
\vdots\\
\dbar \phi_{i,s-1} \\
\dbar \phi_{i,s}\\
\dbar \phi_{i,s+1}
\end{array}\right) + 
\left(\begin{array}{c}
 \omega_{12}^{0,1}\phi_{i2} + \omega_{13}^{0,1}\phi_{i3} + \cdots  + \omega_{1,s+1}^{0,1}\phi_{i,s+1}\\
 \omega_{23}^{0,1}\phi_{i3} + \cdots + \omega_{2,s+1}^{0,1}\phi_{i,s+1}\\
\vdots \\
\omega_{s-1,s}^{0,1}\phi_{i,s} + \omega_{s-1,s+1}^{0,1}\phi_{i,s+1}\\
\omega_{s,s+1}^{0,1}\phi_{i,s+1}\\
 0
\end{array}\right) = 0,
\end{equation}
where $\omega_{ij}^{0,1}$ are the entries of $A^{0,1}$, if we split $A=A^{1,0}+A^{0,1}$ into the parts of 
type (1,0) and (0,1).

This implies that $\phi_{i,s+1}$ is a globally defined (vector-valued) holomorphic function and therefore 
constant for all $i$.
Let us assume that we already know that $\phi_{i,j}$ is constant for all $i$ and for all $j=s+1, s, \ldots, j_0$.
Denote the entry in the $(j_0-1)$-st line by 
$\tau_{j_0-1}^{0,1}:=\omega_{j_0-1,j_0}^{0,1}\phi_{i,j_0}+ \omega_{j_0 -1, j_0+1}^{0,1}\phi_{i,j_0+1}
+ \cdots + \omega_{j_0-1, s+1}^{0,1}\phi_{i,s+1}$.
By (\ref{etaconditions2}) we know that $\partial(\tau_{j_0-1}^{0,1})=0$.
Together with $\dbar(\tau_{j_0-1}^{0,1})=\dbar^2(\phi_{i,j_0-1})=0$ this implies that 
$\tau_{j_0-1}^{0,1}$ is harmonic and in addition $\dbar$-exact. By the Hodge decomposition
theorem for the $\dbar$-operator, saying that 
Ker $(\dbar: \Gamma(X,\mathcal{A}^{0,1}_X) \to \Gamma(X,\mathcal{A}^{0,2}_X))=
\dbar(\Gamma(X,\mathcal{A}^0_X))\oplus {\rm Harm}^{0,1}(X)$, we obtain that $\dbar\phi_{i,j_0-1}=\tau_{j_0-1}^{0,1}=0$,
hence $\phi_{i,j_0-1}$ is constant for all $i$.

Therefore we obtain by induction that all $\phi_{ij}$ are constant. 
Since the $r\times r$-matrix $(\tilde{\phi}_{ij})$
is regular, the differential equation (\ref{diffeq}) implies that $A^{0,1}=0$. Therefore $A$, and hence $\eta$,
is of type (1,0). 
\end{proof}

By Chen's formula the monodromy representation of $(V,\nabla)$ is 
given by 
\begin{equation}
\rho=E_r +\int A + \int A A + \cdots + \int A^s.
\end{equation}

It will suffice to prove the factorization properties of the monodromy representation using the matrix representation
of the monodromy that we obtain with respect to the $D$-flat frame $f$ of the bundle which was chosen at the beginning of
this section.

\begin{lem}  \label{splittinglem}
Let $ A$ be a matrix of 1-forms with a block structure as in (\ref{connectionmatrix}), satisfying $d A=0$ and $A\wedge A=0$.
If $\rho:= E_r + \int A +\int A^2 + \cdots + \int A^{s}$ is a homotopy functional, when restricted to loops based
at $x$, then
every summand $\int A^l$, $l=1,\ldots, s$, has this property, too.
\end{lem}
\begin{proof}
We prove the assertion by induction over $s$. The case $s=1$ is trivial. For $s>1$ we only have to
look at the $(1,s+1)$-entry of $\rho$, since the other entries split into homotopy invariant summands by induction
hypothesis.

The (matrix-valued) entry of $\rho$ in the $(1,s+1)$-place is $\sum_{t=0}^{s-1} I_t$,
where $I_0:=\int\omega_{1,s+1}$ and 
\[I_t:=\sum_{1<i_1 < i_2 < \cdots < i_t < s+1} \int \omega_{1 i_1}\omega_{i_1 i_2}\cdots \omega_{i_t, s+1},
\, \, \, t=1, \ldots, s-1.\]
We want show that every $I_t$ is homotopy invariant when restricted to loops based at $x$.
By Section \ref{MHSsection} we have to show that 
\[(d_C+d_I)(\sum_{1<i_1 < i_2 < \cdots < i_t < s+1} [\omega_{1 i_1}|\omega_{i_1 i_2}|\cdots |\omega_{i_t, s+1}])\in R.\]
But since all $\omega_{ij}$ are closed and because of $A\wedge A=0$, we even obtain
\[
d_I(\sum_{1<i_1 < \cdots < i_t < s+1} [\omega_{1 i_1}|\cdots |\omega_{i_t, s+1}]) =0,
\]
and 
\[d_C(\sum_{1<i_1 <  \cdots < i_t < s+1} [\omega_{1 i_1}|\cdots |\omega_{i_t, s+1}] ) 
=
- \sum_{1<i_1 <  \cdots < i_t < s+1} \sum_{k=1}^{t} \Lambda(i_1,\ldots,i_t, k)=0,\]
where
 \[\Lambda(i_1,\dots,i_t,k):= [\omega_{1 i_1 }|\cdots |\omega_{i_{k-2},i_{k-1}}|\omega_{i_{k-1},i_k}\wedge\omega_{i_k,i_{k+1}} |
 \omega_{i_{k+1},i_{k+2}}|\cdots |\omega_{i_t, s+1}]
\]
with $i_0:=1$ and $i_{t+1}:=s+1$.
\end{proof}

Recall that we defined the ideals
\begin{equation}
I:=J\cap \f^0 + J^2\cap \f^{-1} + \cdots + J^s\cap \f^{-(s-1)}
\end{equation}
and 
\begin{equation}
\bar{I}:=J\cap \bar{\f}^0 + J^2\cap \bar{\f}^{-1} + \cdots + J^s\cap \bar{\f}^{-(s-1)}
\end{equation}
in the introduction, where $\f^i$ denotes the Hodge filtration of $\frac{\C\pi_1(X,x)}{J^{s+1}}$ defined in
(\ref{Hodgefiltration}).

For the proof of the Theorems \ref{theorem} and \ref{bartheorem} stated in the introduction 
we will need the following remark.

\begin{rem} \label{subbundlerem}
{\rm 
As above let $\rho=E_r + \int  A + \int  A^2 + \cdots + \int  A^s$ be the monodromy representation of a flat bundle
factoring through $I$ (or $\bar{I}$) and let $\eta$ be the operator coming from the harmonic category (see
\ref{catexamples}, (5), and \ref{Corlette}), given by the matrix $A$ with respect to the frame chosen above.
By (\ref{etaconditions2}) the matrix $A$ fulfills $\partial  A=0$ and $\dbar A + A\wedge A=0$. 
Denote by $A'$ the matrix belonging
to the $(r_1+\cdots + r_{s-1})$-subbundle and by $A''$ the matrix belonging to the 
$(r_2 + \cdots + r_s)$-quotient bundle. Obviously these satisfy the same formulas as $A$. Therefore they define
the uniquely determined operators on the harmonic bundles belonging to the sub- and quotient bundle in the description
given in Remark \ref{Corlette}. Their monodromy representations
$E_{r-r_{s+1}} + \int A' + \int (A')^2 + \cdots + \int (A')^{s-1}$ and 
$E_{r-r_1} + \int A'' + \int (A'')^2 + \cdots + \int (A'')^{s-1}$
factor through $J^s$ and through $I$ (or $\bar{I}$). 
}
\end{rem}

\noindent {\it Proof of Theorem \ref{theorem}.}

With the notation fixed at the beginning this section we have that 
\begin{equation}
(V,\nabla)= (\oplus_{i=1}^r f_i\mathcal{A}^0_X,  d^{\oplus r} + A )
\end{equation}
with the connection matrix $A$ fulfilling
$\partial A=0$ and $\dbar A +  A\wedge A=0$, and $(\oplus_{i=1}^r f_i\mathcal{A}^0_X, \dbar^{\oplus r} + A)$
being the corresponding Higgs bundle. 

If its underlying holomorphic bundle is trivial, i.e. by Lemma \ref{typelemma} if
$A$ is of type $(1,0)$, the condition $\dbar A +  A\wedge A=0$ implies that 
$\dbar A=0$, since this summand is of type $(1,1)$, whereas the other one is of type $(2,0)$.
Therefore all entries of $A$ are holomorphic 1-forms, hence closed.

The monodromy representation of $(V,\nabla)$ is $\rho=E_r + \int  A + \int  A^2 + \cdots + \int  A^s$.
It follows from the previous lemma that every summand $\int A^l$, $l=1,\ldots, s$, 
is homotopy invariant, when restricted to loops based at $x$. 
Since we have shown that the entries of all matrices $\omega_{ij}$ are holomorphic 1-forms, 
it is clear by (\ref{Hfiltration}) that the integrals of length $l$
are in $\f^l \h^0(B_{s+1}(\mathcal{A}^\bullet))$, hence the factorization property follows immediately.

Conversely, we assume that $\rho=E_r + \int  A + \int  A^2 + \cdots + \int  A^s$ factors through
$I$. By Lemma \ref{typelemma} we have to show that $A$ is necessarily of type $(1,0)$.
We prove this by induction over $s$. 
If $s=0$ the representation $\rho$ is trivial and corresponds to the trivial Higgs bundle 
$(\mathcal{A}^0_X,\dbar)^{\oplus r}$.

We denote by $\eta$ the operator from Remark \ref{Corlette} belonging to $(V,\nabla)$ satisfying 
(\ref{etaconditions}).
Its matrix $A$ with respect to the chosen frame fulfills $\partial  A=0$ and $\dbar A + A\wedge A=0$. 
By Remark \ref{subbundlerem} and induction hypothesis we can assume that all entries of 
$A$ except for those in the $(1,s+1)$-block are of type $(1,0)$.
The condition $\dbar A + A\wedge A=0$ implies that $\dbar\omega_{1,s+1}+(\omega_{12}\wedge\omega_{2,s+1} + \cdots 
+ \omega_{s-1,s}\wedge\omega_{s,s+1}) =0$, from which it follows that $\dbar\omega_{1,s+1}=0$, since the other summands
are of type $(2,0)$, whereas $\dbar\omega_{1,s+1}$ cannot have a $(2,0)$-part. By assumption
it also holds that $\partial\omega_{1,s+1}=0$. From this it follows that $\omega_{1,s+1}$ is harmonic.

In particular $\omega_{1,s+1}$ is closed, hence  
$\int\omega_{1,s+1}$ is a homotopy functional, when restricted to loops based at $x$,  
factoring through $I$ and even factoring through $\f^0(\frac{\C\pi_1(X,x)}{J^{s+1}})$, since integration
over a constant also yields zero.
By the duality of the mixed Hodge structures, $\int\omega_{1,s+1}$ lies in $\f^1\h^0({\rm B}_1(\mathcal{A}^\bullet))$.
Hence there is a (matrix-valued) $\mathcal{C}^\infty$-function $\psi_{1,s+1}$ such that 
$\omega_{1,s+1}-d\psi_{1,s+1}$ is holomorphic, hence also harmonic as $\h^0(X,\Omega^1_X)={\rm Harm}^{1,0}(X)$.
By the Hodge decomposition theorem, saying that 
Ker($d:\Gamma(X,\mathcal{A}^1_X) \to \Gamma(X,\mathcal{A}^2_X))=d(\Gamma(X,\mathcal{A}^0_X))\oplus {\rm Harm}^1(X)$,
we obtain that $d\psi_{1,s+1}$, being harmonic and $d$-exact at the same time, equals zero.
Thus $\omega_{1,s+1}$ is holomorphic, which completes the proof.
\ep

\noindent {\it Proof of Theorem \ref{bartheorem}.}

With the notation above we have that 
\begin{equation}
(V,\nabla)= (\oplus_{i=1}^r f_i\mathcal{A}^0_X,  d^{\oplus r} + A )
\end{equation}
with the connection matrix $A$ fulfilling
$\partial A=0$ and $\dbar A +  A\wedge A=0$, and $(\oplus_{i=1}^r f_i\mathcal{A}^0_X, \dbar^{\oplus r} + A)$
being the corresponding Higgs bundle. 

If its Higgs field is equal to zero, i.e. by Lemma \ref{typelemma}
if $ A$ is of type $(0,1)$, the condition $\partial A=0$ implies that 
all entries of $A$ is anti-holomorphic, hence closed.

The monodromy representation of $(V,\nabla)$ is $\rho=E_r + \int  A + \int  A^2 + \cdots + \int  A^s$.
Since it is homotopy invariant and since $ A$ is closed, it follows from Lemma \ref{splittinglem} that every
$\int A^l$, $l=1,\ldots,s$, is also homotopy invariant, when restricted to loops based at $x$.

Since all 1-forms in the matrix $A$ are anti-holomorphic, it is clear 
that the integrals of length $l$ contain $l$ $d\bar{z}$'s, i.e. 
by definition of the Hodge filtration (\ref{Hfiltration}) they
are in $\bar{\f}^l \h^0(B_{s+1}(\mathcal{A}^\bullet))$, hence the factorization property follows immediately.

Conversely, we assume that $\rho=E_r + \int  A + \int  A^2 + \cdots + \int  A^s$ factors through
$\bar{I}$. By Lemma \ref{typelemma} we have to show that $ A$ is necessarily of type $(0,1)$.
We prove this by induction over $s$.
For $s=0$ the representation $\rho$ is trivial and it corresponds to the trivial Higgs bundle
$(\mathcal{A}^0_X,\dbar)^{\oplus r}$, having Higgs field equal to zero.

By Remark \ref{subbundlerem} and by induction hypothesis we can assume that all entries of $A$ 
except for those in the (matrix-valued)
$(1,s+1)$-block are of type $(0,1)$.
We have to show that $\omega_{1,s+1}$ is harmonic. 

By assumption it holds that $\partial\omega_{1,s+1}=0$. 
It suffices to show that also $\dbar\omega_{1,s+1}=0$. 
We split $\omega_{1,s+1}=\omega_{1,s+1}^{1,0} + \omega_{1,s+1}^{0,1}$ into the
(1,0)- and the (0,1)-part.
From the condition 
$\dbar\omega_{1,s+1}^{1,0} + \dbar\omega_{1,s+1}^{0,1} + (\omega_{12}\wedge\omega_{2,s+1} + \cdots + 
\omega_{1,s}\wedge\omega_{s,s+1})=0$, it follows that
$\dbar\omega_{1,s+1}^{1,0}=0$, since it is the only summand of type (1,1).
Further, $\partial\omega_{1,s+1}^{0,1}=0$ implies that $\omega_{1,s+1}^{0,1}$ is
antiholomorphic, in particular closed. 
Therefore $\dbar\omega_{1,s+1}^{0,1}=d\omega_{1,s+1}^{0,1}-\partial\omega_{1,s+1}^{0,1}=0$.

Thus $\omega_{1,s+1}$ is harmonic, in particular closed, and
we have that $\int\omega_{1,s+1}$ is a homotopy
functional factoring through $\bar{I}$.

Since integration over a constant also yields zero we even know that
$\int_c \omega_{1,s+1}=0$ for all $c \in \bar{\f}^0(\C\pi_1(X,x)/J^2)$, 
i.e. $\int\omega_{1,s+1}\in \bar{\f}^1\h^0({\rm B}_1(\mathcal{A}^\bullet))$ by the duality of the mixed
Hodge structures.

Hence there is a $\mathcal{C}^\infty$-function $\psi_{1,s+1}$ such that 
$\omega_{1,s+1}-d\psi_{1,s+1}$ is anti-holomorphic, i.e. every entry lies in Harm$^{0,1}(X)$.
Since $\omega_{1,s+1}$ is also harmonic, we obtain by the Hodge decomposition theorem,
Ker($d:\Gamma(X,\mathcal{A}^1_X) \to \Gamma(X,\mathcal{A}^2_X))=d(\Gamma(X,\mathcal{A}^0_X))\oplus {\rm Harm}^1(X)$,
that $d\psi_{1,s+1}$, being harmonic and $d$-exact at the time, equals zero.
Therefore $\omega_{1,s+1}$ must be anti-holomorphic, in particular of type (0,1), hence the
operator $\eta$ is of type (0,1).

By Lemma \ref{typelemma} it follows that the Higgs field of the corresponding Higgs bundle $(V,D'' +\eta)$
is equal to zero. \ep

\noindent{\it Acknowledgements.}
I thank M. Nori for suggesting the problem and H. Esnault for the useful discussions we had at different
stages of this work.


\begin{thebibliography}{99}


\bibitem{gm}
{\it W. Goldman, J. Millson}, The deformation theory of representations of
fundamental groups of compact K\"ahler manifolds, Inst. Hautes \'Etudes Sci. Publ. Math.,
{\bf 67} (1988), 43-96

\bibitem{hain0}
{\it R. Hain}, On a generalizaton of Hilbert's 21st problem, Ann. Sci. \'Ecole Norm. Sup. (4),
{\bf 19} (1986), 609-627

\bibitem{hain1}
{\it R. Hain}, The geometry of the mixed Hodge structure on the fundamental group, 
Proc. Symp. Pure Math., {\bf 46} (1987), 247-282

\bibitem{hain2}
{\it R. Hain}, The de Rham homotopy theory of complex algebraic varieties I,
K-Theory, {\bf 1} (1987), 271-324

\bibitem{morgan}
{\it J. Morgan}, The algebraic topology of smooth algebraic varieties, Inst. Hautes \'Etudes
Sci. Publ. Math., {\bf 48} (1978), 137-204

\bibitem{se}
{\it J.-P. Serre}, Lie algebras and Lie groups, Lect. Notes in Math., {\bf 1500}
(1992), Springer

\bibitem{si}
{\it C.T. Simpson}, Higgs bundles and local systems, Inst. Hautes \'Etudes Sci. Publ.
Math., {\bf 75} (1992), 5-95

\end{thebibliography}
\end{document}